\documentclass{ifacconf}

\usepackage{graphics} 
\usepackage{amsmath} \usepackage{amssymb} \usepackage{amsxtra}
\usepackage{arydshln}
\usepackage{amsfonts} \usepackage[mathscr]{eucal}
\usepackage{mathrsfs} \usepackage{setspace} \usepackage{url}
\usepackage[pdftex]{graphicx} \usepackage{multirow}
\usepackage{natbib} 

\def\Real{\mathbb{R}}

      \def\HH{\mathcal{H}}

\newcommand{\Ltwo}{\boldsymbol{\rm L}_{2}}
\newcommand{\Ltwoe}{\boldsymbol{\rm L}_{2e}}

\theoremstyle{definition}

\begin{document}

\begin{frontmatter}

  \title{Converse passivity theorems\thanksref{footnoteinfo}}

  \thanks[footnoteinfo]{The authors gratefully acknowledge the support of the Institute for Mathematics and its Applications, where this work was
    initiated during the 2015-2016 program on Control Theory and its Applications.}

  \author[UMN]{Sei Zhen Khong}    \author[Gro]{Arjan van der Schaft} 

  \address[UMN]{ Institute for Mathematics and its Applications, The University of Minnesota, Minneapolis, MN 55455, USA. (e-mail:
    \texttt{szkhong@umn.edu})}

  \address[Gro]{Johann Bernoulli Institute for Mathematics and Computer Science, University of Groningen, Groningen 9747 AG, Netherlands. (email: \texttt{a.j.van.der.schaft@rug.nl}).}


\begin{abstract}
  Passivity is an imperative concept and a widely utilized tool in the analysis and control of interconnected systems. It naturally arises in the
  modelling of physical systems involving passive elements and dynamics. While many theorems on passivity are known in the theory of robust control,
  very few converse passivity results exist. This paper establishes various versions of converse passivity theorems for nonlinear feedback systems. In
  particular, open-loop passivity is shown to be necessary to ensure closed-loop passivity from an input-output perspective. Moreover, the
  stability of the feedback interconnection of a specific system with an arbitrary passive system is shown to imply passivity of the system itself.
\end{abstract}

\end{frontmatter}

\section{Introduction}
Passivity has emerged to be a crucial concept and tool in the analysis and control of feedback systems, and interconnected systems in general; see
e.g. \cite{Wil72, MoyHil78, Vid81, MegRan97, Sch16}. A salient result is the {\it passivity theorem}, stating (in various versions) that the standard feedback
interconnection of two passive systems is again passive (and hence stable in a certain sense). While passivity theory is deep-rooted in physical
systems modeling and synthesis (such as electrical network theory) based on the essential notions of power and energy, its underlying concepts and
results have turned out to be equally valuable in the broad field of control, ranging from adaptive control to stabilization of nonlinear systems.

While the passivity theorem pervades large parts of systems and control theory, the {\it converse} versions of the passivity theorem seem to be much
less recognized and appreciated. The simplest version of a converse passivity theorem, stating that the feedback interconnection of two systems is
passive {\it if and only if} the two (open-loop) systems are passive, was previously noted and proved within the state-space context in \cite{kerber-vds11}
(see also \cite{Sch16}), and an easy proof in the nonlinear input-output map setting will be provided in Section 3.

The main part of the paper (Section 4) is concerned with a different, and more involved, converse passivity theorem, stating that if the feedback
interconnection of a system with an {\it arbitrary} passive system is {\it stable} (to be specified later on), then the system is necessarily
passive. This basic idea is, sometimes implicitly, underlying a large part of literature on robotics and impedance control; see
e.g. \cite{Str15}. In fact, a version of this result was proved for linear single-input single-output systems in \cite{colgate-hogan} using arguments from
Nyquist stability theory, exactly with this motivation in mind. In robotics the motivation for this converse passivity theorem can be formulated as
follows. Consider a controlled robot interacting in operation with its environment (the normal scenario). In many cases the environment is largely
unknown, while at the same time the {\it stability} of the robot interacting with its environment can be often considered as a {\it sine qua
  non}. Since the interaction of the robot with its environment typically takes place via the conjugated variables of (generalized) velocity and
force, whose product is equal to {\it power}, it is not completely unreasonable to assume that the environment, seen from the interaction port with
the robot, is, although unknown, a {\it passive system}. Thus the converse passivity theorem treated in Section 4 gives a clear rationale for the
often expressed design and control principle \citep{Str15, colgate-hogan} that a controlled robot should be passive at its interaction port with the
unknown environment\footnote{The second author would like to thank Stefano Stramigioli for many inspiring conversations on this and related
  subjects.}. Differently from \cite{colgate-hogan, Str15}, the proof of the general nonlinear converse passivity theorem treated in Section 4 will be based on
the S-procedure lossless theorem (see \cite{MegTre93} or \cite[Thm. 7]{Jon01}). 

\section{Notation and preliminaries}

Denote by $\Ltwo^n$ the set of $\Real^n$-valued Lebesgue square-integrable functions:
\begin{align*}
\Ltwo^n := \Big\{v : [0, \infty) & \to \Real^n : \\
& \|v\|_2^2 = \langle v, v \rangle := \int_0^\infty v(t)^Tv(t) \, dt < \infty\Big\}.
\end{align*}
Define the truncation operator
\[
(P_Tv)(t) := \left\{\begin{array}{ll}
v(t) & \quad t \in [0, T) \\
0 & \quad \text{otherwise},
\end{array} \right.
\]
and the extended spaces
\begin{align*}
 \Ltwoe^n & := \{v : [0, \infty) \to \Real^n : P_Tv \in \Ltwo \; \forall T \in [0, \infty)\}.
\end{align*}
In what follows, the superscript $n$ is often suppressed for notational simplicity.  Define the shift operator $(S_T u)(t) = u(t - T)$ for $T \geq 0$
and denote the identity system by $I : \Ltwoe \to \Ltwoe$.  A system $\Delta : \Ltwoe \to \Ltwoe$ is said to be \emph{causal} if
$P_T \Delta P_T = P_T\Delta$ for all $T \geq 0$. It is said to be \emph{time-invariant} if $S_T \Delta = \Delta S_T$ for all $T > 0$. A causal
$\Delta$ is called \emph{bounded} if its Lipschitz bound~\cite[Section 2.4]{Wil71} is finite, i.e.
\begin{align*} 
\|\Delta\| := \sup_{T > 0; \|P_Tu\|_2 \neq 0} \frac{\|P_T \Delta u\|_2}{\|P_T u\|_2} = \sup_{0 \neq u \in \Ltwo} \frac{\|\Delta u\|_2}{\|u\|_2} < \infty.
\end{align*}
$\Delta$ is said to be \emph{passive}~\citep{Wil72, Sch16} if
\begin{align} \label{eq: passive}
\inf_{T > 0, u \in \Ltwoe} \int_0^T u(t)^T (\Delta u)(t) \, dt \geq 0,
\end{align}
\emph{strictly passive} if there exists $\epsilon > 0$ such that
\[
\int_0^T u(t)^T (\Delta u)(t) \, dt \geq \epsilon (\|P_T u\|_2^2 + \|P_T \Delta u\|_2^2)
\]
$\forall u \in \Ltwoe, T >0$, and \emph{output strictly passive} if there exists an $\epsilon > 0$ such that
\[
\int_0^T u(t)^T (\Delta u)(t) \, dt \geq \epsilon \|P_T \Delta u\|_2^2 \quad \forall u \in \Ltwoe, T >0.
\]

\begin{lem} \label{lem: passive}
If $\Delta$ is bounded, then passivity is equivalent to 
\begin{align} \label{eq: bounded_passive}
\inf_{u \in \Ltwo} \int_0^\infty u(t)^T (\Delta u)(t) \, dt \geq 0.
\end{align}
Similarly, strict passivity and output strict passivity of $\Delta$ are equivalent to
\[
\int_0^\infty u(t)^T (\Delta u)(t) \, dt \geq \epsilon (\|u\|_2^2 + \|\Delta u\|_2^2) \quad \forall u \in \Ltwo
\]
and
\[
\int_0^\infty u(t)^T (\Delta u)(t) \, dt \geq \epsilon \|\Delta u\|_2^2 \quad \forall u \in \Ltwo,
\]
respectively.
\end{lem}

\begin{pf}
  First, note that \eqref{eq: bounded_passive} can be obtained from \eqref{eq: passive} by restricting $u \in \Ltwo$ and taking $T$ to
  $\infty$. Conversely, suppose that \eqref{eq: passive} does not hold, then there exist $T > 0$ and $u \in \Ltwoe$ such that
  $\int_0^T u(t)^T (\Delta u)(t) \, dt < 0$. Let $\bar{u}(t) := u(t)$ for $t \in [0, T)$ and $\bar{u}(t) := 0$ for $t \geq T$. Then
  $\bar{u} \in \Ltwo$ and $\int_0^\infty \bar{u}(t)^T (\Delta \bar{u})(t) = \int_0^T u(t)^T (\Delta u)(t) \, dt < 0$. That is, \eqref{eq:
    bounded_passive} is violated. This completes the proof for the first part of the lemma. The rest of the lemma can be shown in a similar fashion. \hfill$\qed$
\end{pf}

\setlength{\unitlength}{0.7cm}
\begin{figure}[h]
  \centering \includegraphics[scale=0.6]{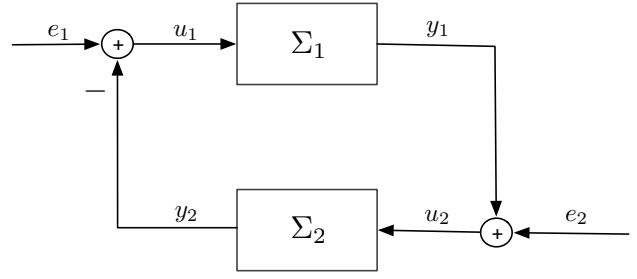}
  \caption{Feedback configuration} \label{fig: feedback}
\end{figure}

The main object of study in this paper is the feedback interconnection of causal systems $\Sigma_1 : \Ltwoe \to \Ltwoe$ and $\Sigma_2 : \Ltwoe \to \Ltwoe$ described by
\begin{align} \label{eq: FB}
\begin{split}
 u_1 & = e_1 - y_2; \quad u_2 = e_2 + y_1; \\
 y_1 & = \Sigma_1u_1; \qquad y_2 = \Sigma_2u_2;
\end{split}
\end{align}
see Figure~\ref{fig: feedback}.

\begin{defn} \label{def: FB} The feedback interconnection $\Sigma_1 \| \Sigma_2$ is said to be \emph{well-posed} if the map $(y_1, y_2) \mapsto (e_1, e_2)$ defined by \eqref{eq: FB}
  has a causal inverse $F$ on $\Ltwoe$. It is \emph{finite-gain stable} if in addition to being well-posed, $\|F\| < \infty$. 
A well-posed $\Sigma_1 \| \Sigma_2$ is said to be \emph{passive} if the map $(e_1, e_2) \mapsto (y_1, y_2)$ is passive. In the case where $e_2 = 0$,
the feedback interconnection is said to be finite-gain stable if it is well-posed and $e_1 \mapsto y_1$ is bounded.
\end{defn}

\section{Feedback passivity}

In this section, a simple proof that passivity of the closed-loop system implies passivity of the open-loop components is provided. Contrary to the
state-space setting in \cite{kerber-vds11}, the result adopts the input-output perspective and is applicable to infinite-dimensional systems, such as
those modelled by time-delay or partial differential equations.

\begin{thm}
  Given causal $\Sigma_1$ and $\Sigma_2$ for which $\Sigma_1 \| \Sigma_2$ is well posed, $\Sigma_1 \| \Sigma_2$ is passive if, and only if,
  $\Sigma_1$ and $\Sigma_2$ are passive.
\end{thm}

\begin{pf}
 ($\Longrightarrow$) By hypothesis, 
\[
\inf_{T > 0, e_1, e_2 \in \Ltwoe} \int_0^T e_1(t)^T y_1(t) + e_2(t)^Ty_2(t) \, dt \geq 0.
\]
Using \eqref{eq: FB} and the well-posedness of $\Sigma_1 \| \Sigma_2$, this is equivalent to
\begin{align*}
  0 & \leq \inf_{T > 0, u_1, u_2 \in \Ltwoe} \int_0^T (u_1(t) + y_2(t))^T y_1(t) + \\
 & \qquad\qquad\qquad\qquad\qquad\qquad (u_2(t) - y_1(t))^Ty_2(t) \, dt \\
    & = \inf_{T > 0, u_1, u_2 \in \Ltwoe} \int_0^T u_1(t)^T y_1(t) + u_2(t)^Ty_2(t) \, dt
\end{align*}
This implies that 
\[
\inf_{T > 0, u_1\in \Ltwoe} \int_0^T u_1(t)^T y_1(t) \, dt \geq 0
\]
(by setting $u_2 = 0$) and 
\[
\inf_{T > 0, u_2\in \Ltwoe} \int_0^T u_2(t)^T y_2(t) \, dt \geq 0
\]
(by setting $u_1 = 0$), which are equivalent to passivity of $\Sigma_1$ and $\Sigma_2$, respectively.

($\Longleftarrow$) This direction is well-known in the literature and can be shown by reversing the arguments above. \hfill$\qed$
\end{pf}

\section{Passivity as a necessary and sufficient condition for stable interaction}

In this section we show that a necessary and sufficient condition in order that the closed-loop system arising from interconnecting a given system to an unknown, but passive, system is stable, is that the system is passive itself. We will formulate three slightly different versions of this main result.

\begin{thm} \label{thm: passive}
  Given a bounded time-invariant $\Sigma_1$, the feedback interconnection $\Sigma_1 \| \Sigma_2$ is finite-gain stable for all bounded passive $\Sigma_2$ if, and only if,
  $\Sigma_1$ is strictly passive.
\end{thm}

\begin{pf}
 Sufficiency is well known in the literature. Indeed, by Lemma~\ref{lem: passive}, strict passivity of $\Sigma_1$ and passivity of $\Sigma_2$ yield
\begin{align*}
\epsilon (\|y_1\|_2^2 + \|u_1\|_2^2) & \leq \langle u_1, y_1 \rangle + \langle u_2, y_2 \rangle \\
& = \langle e_1 - y_2, y_1 \rangle + \langle e_2 + y_1, y_2 \rangle \\
& = \langle e_1, y_1 \rangle + \langle e_2, y_2 \rangle.
\end{align*}
Therefore,
\[
\|y_1\|_2^2 + \langle e_1 - y_2, e_1 - y_2 \rangle \leq \frac{1}{\epsilon} (\langle e_1, y_1 \rangle + \langle e_2, y_2 \rangle)
\]
or
\[
\|y_1\|_2^2 + \|y_2\|_2^2 -2 \langle e_1, y_2 \rangle + \|e_1\|_2^2 \leq \frac{1}{\epsilon} (\langle e_1, y_1 \rangle + \langle e_2, y_2 \rangle)
\]
It follows that
\[
\|y\|_2^2 \leq 2 \langle e_1, y_2 \rangle + \frac{1}{\epsilon} \langle e, y \rangle \leq \left(2 + \frac{1}{\epsilon}\right) \|e\|_2\|y\|_2,
\]
where $y := (y_1, y_2)^T$ and $e := (e_1, e_2)^T$ and the Cauchy-Schwarz inequality has been used.

To show necessity, define 
\[
\HH := \{h = (u_1, u_2, e_1, e_2) \in \Ltwo\ |\ u_2 = e_2 + \Sigma_1 u_1\}.
\]
Note that if $h \in \HH$, then $S_T h \in \HH$ for all $T \geq 0$ due to the time-invariance of $\Sigma_1$. Define the quadratic forms $\sigma_i : \HH
\to \Real$, $i = 0, 1$, as 
\begin{align*}
 \sigma_0(u_1, u_2, e_1, e_2) & := \left\langle 
\begin{bmatrix} 
 u_1 \\
 u_2 \\
 e_1 \\
 e_2
\end{bmatrix},
\begin{bmatrix}
I & 0 & 0 & 0 \\
0 & I & 0 & 0 \\
0 & 0 & -\gamma I & 0 \\
0 & 0 & 0 & -\gamma I
\end{bmatrix}
\begin{bmatrix} 
 u_1 \\
 u_2 \\
 e_1 \\
 e_2
\end{bmatrix}
\right\rangle \\
 \sigma_1(u_1, u_2, e_1, e_2) & := \left\langle 
\begin{bmatrix} 
 u_1 \\
 u_2 \\
 e_1 \\
 e_2
\end{bmatrix},
\frac{1}{2} \begin{bmatrix}
0 & -I & 0 & 0 \\
-I & 0 & I & 0 \\
0 & I &0 & 0 \\
0 & 0 & 0 & 0
\end{bmatrix}
\begin{bmatrix} 
 u_1 \\
 u_2 \\
 e_1 \\
 e_2
\end{bmatrix}
\right\rangle.
\end{align*}
By Lemma~\ref{lem: passive}, stability of
$\Sigma_1 \| \Sigma_2$ for all bounded passive $\Sigma_2$ implies the existence of $\gamma >0$ such that
\begin{align*}
 \sigma_0(u_1, u_2, e_1, e_2) \leq 0 \quad \forall & (u_1, u_2, e_1, e_2) \in \HH \\
 & \text{such that} \quad \sigma_1(u_1, u_2, e_1, e_2) \geq 0.
\end{align*}
This is equivalent, via the S-procedure lossless theorem (see \cite{MegTre93} or \cite[Thm. 7]{Jon01}), to the existence of $\tau \geq 0$ such that
\begin{align*}
 \sigma_0(u_1, u_2, e_1, e_2) + \tau \sigma_1(u_1, u_2, & e_1, e_2) \leq 0 \\
 & \forall (u_1, u_2, e_1, e_2) \in \HH.
\end{align*}
In the subset $\{(u_1, u_2, 0, 0) \in \Ltwo\ |\ u_2 = \Sigma_1 u_1\} \subset \HH$, this implies that
\[
\|\Sigma_1 u_1\|_2^2 + \|u_1\|_2^2 - \tau \langle u_1, \Sigma_1 u_1 \rangle \leq 0 \quad \forall u_1 \in \Ltwo.
\]
Equivalently,
\[
\tau \langle u_1, \Sigma_1 u_1 \rangle \geq \|\Sigma_1 u_1\|_2^2 + \|u_1\|_2^2 \quad \forall u_1 \in \Ltwo.
\]
It is obvious from the inequality above that $\tau \neq 0$, hence
\[
\langle u_1, \Sigma_1 u_1 \rangle \geq \frac{1}{\tau} (\|\Sigma_1 u_1\|_2^2 + \|u_1\|_2^2) \quad \forall u_1 \in \Ltwo,
\]
i.e. $\Sigma_1$ is strictly passive, by Lemma~\ref{lem: passive}. \hfill$\qed$
\end{pf}

Certainly from a state space point of view the above theorem has the drawback of relying on {\it strict} passivity, since it is known that input
strict passivity can only occur for systems with direct feedthrough terms. Thus, it excludes a large class of physical systems. The following
alternative version avoids this problem by relying only on {\it output} strict passivity.

\setlength{\unitlength}{0.7cm}
\begin{figure}[h]
  \centering \includegraphics[scale=0.55]{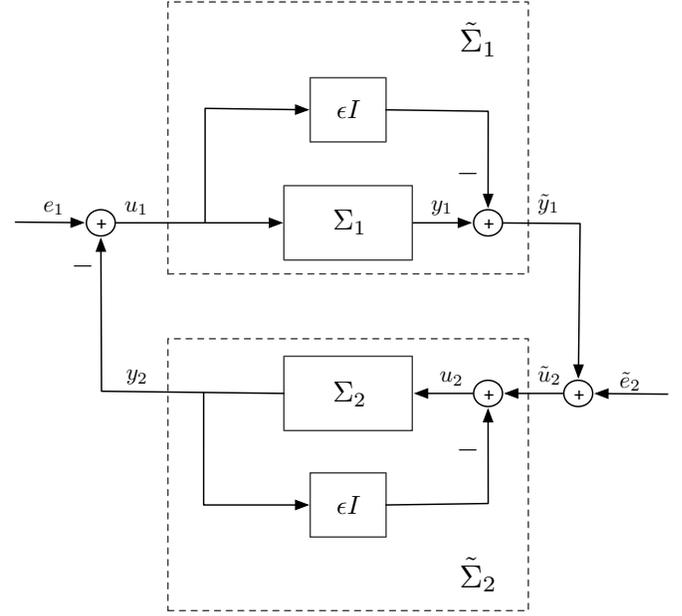}
  \caption{Loop transformation} \label{fig: loop_trans}
\end{figure}

\begin{thm}
  Given a bounded time-invariant $\tilde{\Sigma}_1$, the feedback interconnection $\tilde{\Sigma}_1 \| \tilde{\Sigma}_2$ is finite-gain stable for all
  output strictly passive $\tilde{\Sigma}_2$ if, and only if, $\tilde{\Sigma}_1$ is output strictly passive.
\end{thm}

\begin{pf}
  Sufficiency is well known in the literature and can be shown in a similar manner using the arguments in the sufficiency proof for Theorem~\ref{thm:
    passive}. For necessity, note that any output strictly passive $\tilde{\Sigma}_2$ can be written as the negative feedback interconnection of
  $\Sigma_2$ and $\epsilon I$ for some bounded passive $\Sigma_2$ and $\epsilon > 0$; see Figure~\ref{fig: loop_trans}. To see this, let $\Sigma_2$ be
  bounded and observe by Lemma~\ref{lem: passive} that $\tilde{\Sigma}_2$ satisfies
\[
\int_0^\infty \tilde{u}_2(t)^T y_2(t) \, dt \geq \epsilon \|y_2\|_2^2 \quad \forall \tilde{u}_2 \in \Ltwo,
\]
whereby
\[
\int_0^\infty (\tilde{u}_2(t) - \epsilon y_2(t))^T y_2(t) \, dt = \int_0^\infty u_2(t)^T y_2(t) \, dt \geq 0.
\]
The last inequality holds for all $u_2 \in \Ltwo$, since given any $u_2 \in \Ltwo$, $\tilde{u}_2 := (I + \epsilon \Sigma_2)u_2 \in \Ltwo$ yields the
desired $u_2$. Therefore, it follows that $\Sigma_2$ is passive. By the same token, the negative feedback interconnection of a bounded passive
$\Sigma_2$ and $\epsilon I$ with $\epsilon > 0$ is necessarily output strictly passive. 

By defining $\Sigma_1 := \tilde{\Sigma}_1 + \epsilon I$ as illustrated in Figure~\ref{fig: loop_trans}, one obtains the loop transformation
configuration therein. Consequently, the finite-gain stability of the feedback interconnection $\tilde{\Sigma}_1 \| \tilde{\Sigma}_2$ in
Figure~\ref{fig: loop_trans} is equivalent to that of $\Sigma_1 \| \Sigma_2$ in Figure~\ref{fig: feedback}~\citep[Section 3.5]{GreLim95}. Application
of Theorem~\ref{thm: passive} then yields that $\Sigma_1$ is strictly passive. For sufficiently small $\epsilon > 0$, it then follows that
$\tilde{\Sigma}_1 = \Sigma_1 - \epsilon I$ is output strictly passive. \hfill$\qed$
\end{pf}

Another feature of Theorem \ref{thm: passive} is the fact that it requires an external signal $e_2$, which is not the typical case in robotics
applications. This motivates the following version of converse passivity theorem. Recall that an output strictly passive system has finite $L_2$-gain
\citep{Sch16}.

\begin{thm}
Given a bounded time-invariant $\Sigma_1$, the feedback interconnection $\Sigma_1 \| \Sigma_2$ (with $e_2=0$) has finite $L_2$-gain from $e_1$ to $y_1$ for all passive $\Sigma_2$ if, and only if, $\Sigma_1$ is output strictly passive.
\end{thm}

\begin{pf}
 Sufficiency is well known in the literature. Indeed, if $\Sigma_1$ is output strictly passive and $\Sigma_2$ is passive, then for some $\varepsilon > 0$
 \begin{align*}
 \langle e_1, y_1 \rangle & = \langle u_1 + y_2, y_1 \rangle \\
 & = \langle u_1, y_1 \rangle + \langle y_2, y_1 \rangle \\
 & = \langle u_1, y_1 \rangle + \langle u_2, y_2 \rangle \geq \varepsilon \| y_1\|_2^2,
 \end{align*}
 showing that the closed-loop system is $\varepsilon$-output strictly passive, and hence has $L_2$-gain $\leq \frac{1}{\varepsilon}$. 

 To show necessity, define
\[
\HH := \{h = (u_1, y_1, e_1) \in \Ltwo\ |\ y_1 = \Sigma_1 u_1\}.
\]
Note that if $h \in \HH$, then $S_T h \in \HH$ for all $T \geq 0$ due to the time-invariance of $\Sigma_2$. Define the quadratic forms $\sigma_i : \HH
\to \Real$, $i = 0, 1$, as 
\begin{align*}
 \sigma_0(u_1, y_1, e_1) & := \frac{1}{2}\left\langle 
\begin{bmatrix} 
 u_1 \\
 y_1 \\
 e_1 
\end{bmatrix},
\begin{bmatrix}
0 & 0 & 0  \\
0 & I & 0  \\
0 & 0 & -\gamma^2 I 
\end{bmatrix}
\begin{bmatrix} 
 u_1 \\
 y_1 \\
 e_1 
\end{bmatrix}
\right\rangle \\
 \sigma_1(u_1, y_1, e_1) & := \frac{1}{2} \left\langle 
\begin{bmatrix} 
 u_1 \\
 y_1 \\
 e_1 
\end{bmatrix},
\begin{bmatrix}
0 & -I & 0 \\
-I & 0 & I \\
0 & I &0  
\end{bmatrix}
\begin{bmatrix} 
 u_1 \\
 y_1 \\
 e_1 
\end{bmatrix}
\right\rangle.
\end{align*}
By Lemma~\ref{lem: passive}, stability of
$\Sigma_1 \| \Sigma_2$ for all bounded passive $\Sigma_1$ implies the existence of $\gamma $ such that
\begin{align*}
 \sigma_0(u_1, y_1, e_1) \leq 0 \quad \forall & (u_1, y_1, e_1) \in \HH \\
 & \text{such that} \quad \sigma_1(u_1, y_1, e_1) \geq 0.
\end{align*}
This is equivalent, via the S-procedure lossless theorem (see \cite{MegTre93} or \cite[Thm. 7]{Jon01}), to the existence of $\tau \geq 0$ such that
\[
 \sigma_0(u_1, y_1, e_1) + \tau \sigma_1(u_1, y_1, e_1) \leq 0 \quad \forall (u_1, y_1, e_1) \in \HH.
\]
This implies that
\[
- \tau \langle u_1, y_1 \rangle+  \tau \langle e_1, y_1 \rangle+ \frac{1}{2} \| y_1 \|_2^2 - \frac{1}{2}\gamma^2  \|e_1\|_2^2  \leq 0 \quad \forall e_1 \in \Ltwo ,
\]
and thus in the subset $\{(u_1, y_1, 0) \in \Ltwo\ |\ y_1 = -\Sigma_1 u_1 \} \subset \HH$, this yields
\[
\langle u_1, y_1 \rangle \geq \frac{1}{2\tau} \| y_1 \|_2^2 \quad \forall u_1 \in \Ltwo,
\]
i.e., $\Sigma_2$ is output strictly passive, by Lemma~\ref{lem: passive}. \hfill$\qed$
\end{pf}
In the case of {\it linear single-input single-output} systems, a version of the above theorem was proved before in \cite{colgate-hogan}, using an
argument based on Nyquist stability theory\footnote{Roughly speaking, by showing that if $\Sigma_1$ is {\it not} passive, a positive-real transfer
  function (corresponding to a passive system $\Sigma_2$) can be constructed such that the closed-loop system fails the Nyquist stability test;
  cf. \cite{colgate-hogan}.}.

\section{Conclusions}

Several versions of converse passivity theorems for nonlinear systems are provided. Besides contributing to robust control theory, these fundamental
results have implications in the field of robotics, as described in the introduction. Future work will involve seeking similar results within the
context of large-scale interconnected systems.

\end{document}